\documentclass[11pt]{amsart}
\usepackage[mathscr]{eucal}
\usepackage{amsfonts}
\usepackage{amsmath}
\usepackage{amsthm}
\usepackage{amssymb}
\usepackage{amscd}
\usepackage{latexsym}

\theoremstyle{plain}
\newtheorem{definition}[equation]{Definition}

\newtheorem{corollary}[equation]{Corollary}

\newtheorem{proposition}[equation]{Proposition}
\newtheorem{theorem}[equation]{Theorem}

\theoremstyle{definition}

\numberwithin{equation}{subsection}

\begin{document}

\title{Hochschild and cyclic homology of central extensions of preprojective algebras of ADE quivers}
\author{Ching-Hwa Eu}
\maketitle
\pagestyle{myheadings}
\markboth{Ching-Hwa Eu}{Hochschild and cyclic homology of central extensions of preprojective algebras ...}
\section{Introduction}
Preprojective algebras of quivers were introduced by Gelfand and
Ponomarev in 1979 in order to provide a model for quiver representations
(in the special case of finite Dynkin quivers). Since then,
preprojective algebras have found many other important applications, see
e.g. \cite{CBH}. Ironically, it is exactly in the case of finite
Dynkin quivers, originally considered by Gelfand and Ponomarev,
that preprojective algebras fail to have certain good properties
enjoyed by the preprojective algebras of other connected quivers;
for instance, their deformed versions are not flat. Motivated by
this, the paper \cite{ER} introduces central extensions
of preprojective algebras of finite Dynkin quivers,
and shows that they have better properties,
in particular their deformed versions are flat. The following
paper \cite{ELR} computes the center $Z$ and the trace space
$A/[A,A]$ for the deformed preprojective algebra $A$; the answer
turns out to be related to the structure of the maximal nilpotent
subalgebra of the simple Lie algebra attached to the quiver.

The goal of this paper is to generalize the results of \cite{ELR}
by calculating the additive structure of the Hochschild homology and cohomology
of
$A$ and the cyclic homology of $A$, and to describe the universal
deformation of $A$. Namely, we show that the (co)homology is
periodic with period 4, and compute the first four (co)homology groups in each
case. We plan to study the product structures on the (co)homology
in a separate publication.

We note that the Hochschild cohomology of usual preprojective
algebras $A_0$ (without central extension), together with the cup
product, is studied in \cite{ES} (in the case of type A).
This is done by using the periodic resolution of $A_0$
with period 6 constructed by Schofield. This leads to the
cohomology being periodic with period 6. We note that our
periodic resolution with period 4 for $A$ is quite similar
to the Schofield resolution.\footnote{The reason our resolution has smaller
period than Schofield's is that $A$, unlike $A_0$, has a
symmetric invariant pairing, and hence the Nakayama automorphism
of $A$, unlike that of $A_0$, is the identity.}
In a separate paper, we will apply our methods to
computing the Hochschild (co)homology and cyclic homology of preprojective algebras of type D and E.

The structure of the paper is as follows.
In Section 2 we discuss preliminaries.
In Section 3, we define the periodic
resolution with period 4 for A, and use it to compute the
Hochschild homology and cohomology of $A$. In Section 4,
we compute the cyclic homology of A, and find the Hilbert series
for all the homology and cohomology, using the results
of \cite{ELR} and combinatorial identities from \cite{RS}.
Finally, in Section 5 we use the result about $HH^2(A)$ to find a
universal deformation of $A$ (the deformation theory of $A$ turns
out to be unobstructed).

{\bf Acknowledgments.} The author would like to thank his adviser
Pavel Etingof for posing the problem and many useful discussions
and explanations.
\section{Preliminaries}
\subsection{Quivers, path algebras and preprojective algebras}
Let $Q$ be a quiver of ADE type with vertex set $I$ and $|I|=r$. We define $Q^*$ to be the quiver obtained from $Q$ by reversing all of its arrows. We call $\bar Q=Q\cup Q^*$ the \emph{double} of $Q$.

The concatenation of these arrows generate the \emph{nontrivial paths} inside the quiver $\bar Q$. We define $e_i$, $i\in I$ to be the \emph{trivial path} which starts and ends at $i$. The \emph{path algebra} $P=\mathbb{C}\bar Q$ of $Q$ over $\mathbb{C}$ is the $\mathbb{C}$-algebra with basis the paths in $\bar Q$ and the product $xy$ of two paths $x$ and $y$ to be their concatenation if they are compatible and $0$ if not. We define the \emph{Lie bracket} $[x,y]=xy-yx$. By taking the quotient $P/(\sum\limits_{a\in Q}[a,a^*]),$ we obtain the \emph{preprojective algebra} denoted by $\Pi_Q$.

Let $R=\bigoplus\limits_{i\in I}\mathbb{C}e_i$. Then $P$ (and therefore $\Pi_Q$) is naturally an $R$-bimodule.

\subsection{The symmetric bilinear form, roots and weights}
We write $a\in Q$ to say that $a$ is an arrow in $Q$. Let $h(a)$ denote its  \emph{head} and $t(a)$ its \emph{tail}, i.e. for $a:i\rightarrow j$, $h(a)=j$ and $t(a)=i$. The \emph{Ringel form} of $Q$ is the bilinear form on $\mathbb{Z}^I$ defined by
\[
\langle\alpha,\beta\rangle=\sum\limits_{i\in I}\alpha_i\beta_i-\sum\limits_{a\in Q}\alpha_{t(a)}\beta_{h(a)}
\]
for $\alpha,\,\beta\in\mathbb{Z}^I$. We define the \emph{quadratic form} $q(\alpha)=\langle\alpha,\alpha\rangle$ and the \emph{symmetric bilinear form}  $(\alpha,\beta)=\langle\alpha,\beta\rangle+\langle\beta,\alpha\rangle$.
It can be shown that $q$ is positive definite for a finite Dynkin quiver $Q$.

We define the \emph{set of roots} $\Delta=\{\alpha\in\mathbb{Z}^I|q(\alpha)=1\}$.

We call the elements of $\mathbb{C}^I$ \emph{weights}. A weight $\mu=(\mu_i)$ is called \emph{regular} if the inner product $(\mu,\alpha)\neq 0$ for all $\alpha\in\Delta$. We call the coordinate vectors $\varepsilon_i\in\mathbb{C}^I$ the \emph{fundamental weights} and define $\rho$ to be the sum of all fundamental weights.

We call $h=\frac{|\Delta|}{|I|}$ the \emph{Coxeter number} of $Q$.
\subsection{The centrally extended preprojective algebra}
Let $\mu=(\mu_i)$ be a regular weight. 
We define the \emph{centrally extended preprojective algebra} $A=A^\mu$ to be the quotient of $P[z]$ ($z$ is a central variable) by the relation $\sum\limits_{a\in Q}[a,a^*]=z(\sum\limits_{i\in I}\mu_ie_i)$. By taking the quotient $A/(z)$, we obtain the usual \emph{preprojective algebra} $\Pi_Q=P/(\sum\limits_{a\in Q}[a,a^*])$.

The grading on $A$ is given by $\deg(R)=0$, $\deg(a)=\deg(a^*)=1$ and $\deg(z)=2$.

From now on, we assume $\mu$ to be a generic weight or $\mu=\rho$.

\subsection{The Hilbert series}
\begin{definition} \textnormal{(The Hilbert series of vector spaces)}\\
Let $W=\bigoplus_{d\geq0}W(d)$ be a $\mathbb Z_+$-graded
vector space, with finite dimensional homogeneous subspaces. 
We define the \emph{Hilbert series} $h_W(t)$ to be the series
\begin{displaymath}
h_W(t)=\sum\limits_{d=0}^{\infty}\dim W[d]t^d.
\end{displaymath}
\end{definition}
\begin{definition} \textnormal{(The Hilbert series of bimodules)}\\
Let $W=\bigoplus_{d\geq0}W(d)$ be a $\mathbb{Z_+}$-graded bimodule over the ring $R=\bigoplus_{i\in I}\mathbb{C}$ ($I$ is a finite set), so we can write $W=\bigoplus_{i,j\in I} W_{i,j}$. We define the \emph{Hilbert series} $H_W(t)$ to be a matrix valued series with the entries 
\begin{displaymath}
H_W(t)_{i,j}=\sum\limits_{d=0}^{\infty}\dim\ W(d)_{i,j}t^d.
\end{displaymath}
\end{definition}
\section{Hochschild homology/cohomology and cyclic homology of A}
\subsection{Periodic projective resolution of A}
Let $V$ be the $R$-bimodule which is generated by the arrows in $\bar Q$ (i.e. the degree $1$-part of $A$). For a $\mathbb{Z}-$graded $R$-bimodule $M$, we denote $M[i]$ to be the bimodule $M$, shifted by degree $i$ (i.e. $M(d)=M[i](d+i)$).\\

We want to compute Hochschild homology and cohomology of $A$, so we want to find a projective resolution of $A$.

Let 
\begin{align*}
C_{-1}&=A,\\
C_0&=A\otimes_RA,\\
C_1&=(A\otimes_RV\otimes_RA)\oplus(A\otimes_RA)[2],\\ C_2&=(A\otimes_RV\otimes_RA)[2]\oplus(A\otimes_RA)[2],\\
C_3&=A\otimes_RA[4],\\
C_4&=C_0[2h].
\end{align*}
We define the following $A-$bimodule-homomorphisms $d_i:C_i\rightarrow C_{i-1}$:\\
\[d_0(b_1\otimes b_2)=b_1b_2,\]
\[d_1(b_1\otimes\alpha\otimes b_2,b_3\otimes b_4)=b_1\alpha\otimes b_2-b_1\otimes\alpha b_2+b_3z\otimes b_4-b_3\otimes zb_4,\]
\begin{multline*}
d_2(b_1\otimes\alpha\otimes b_2,b_3\otimes b_4)=
(-b_1z\otimes\alpha\otimes b_2+b_1\otimes\alpha\otimes zb_2+\sum\limits_{a\in\bar Q}\epsilon_ab_3a\otimes a^*\otimes b_4\\
+\sum\limits_{a\in\bar Q}\epsilon_ab_3\otimes a\otimes a^*b_4,-b_3\mu\otimes b_4+b_1\alpha\otimes b_2-b_1\otimes \alpha b_2),
\end{multline*}
where we introduce the notation $\epsilon_a=\left\{\begin{array}{cc}+1&a\in Q\\-1&a\in Q^*
\end{array},
\right.$
\[
d_3(b_1\otimes b_2)=(\sum\limits_{a\in\bar Q}\epsilon_ab_1a\otimes a^*\otimes b_2+\sum\limits_{a\in\bar Q}\epsilon_ab_1\otimes a\otimes a^*b_2,b_1z\otimes b_2-b_1\otimes zb_2),
\]
\[
d_4(b_1\otimes b_2)=\sum b_1x_i\otimes x_i^*b_2,
\] where $\{x_i\}$ is a basis of $A$ and $\{x_i^*\}$ the dual basis under the (symmetric and nondegenerate) trace form $(x,y)=Tr(xy)$ introduced in \cite[Section 2.2.]{ELR}. It is easy to see that $d_4$ is independent of the choice of the basis $\{x_i\}$. It is clear that all $d_i$ are degree-preserving.

Using the trace form, it is easy to show that $\sum ax_i\otimes x_i^*=\sum x_i\otimes x_i^*a$ for any $a\in A$:
\[\sum ax_i\otimes x_i^*=\sum\sum (ax_i,x_j^*)x_j\otimes x_i^*=
\sum\sum x_i\otimes (x_i^*a,x_j)x_j^*=\sum x_i\otimes x_i^*a.\]
This implies 
\[
 d_4(b_1\otimes b_2)=b_0(b_1\otimes b_2)\sum x_i\otimes x_i^*
\]

\begin{theorem}
From the maps $d_i$ we obtain the following projective resolution $C_\bullet$ of $A$ with period $4$:
\[
\cdots\stackrel{d_3[2h]}{\rightarrow}C_2[2h]\stackrel{d_2[2h]}{\rightarrow}C_1[2h]\stackrel{d_1[2h]}{\rightarrow}C_0[2h]\stackrel{d_4}{\rightarrow} C_3\stackrel{d_3}{\rightarrow}C_2\stackrel{d_2}{\rightarrow}C_1\stackrel{d_1}{\rightarrow}C_0\stackrel{d_0}{\rightarrow}A\rightarrow 0.
\]
\end{theorem}
\begin{proof}
Let us first show that these $C_i,\,d_i$ define a complex. We show that $d_id_{i+1}=0$ for $i\leq 3$ and $d_4d_1[2h]=0$:
\[d_0d_1(b_1\otimes\alpha\otimes b_2, b_3\otimes b_4)=d_0(b_1\alpha\otimes b_2-b_1\otimes\alpha b_2+b_3z\otimes b_4-b_3\otimes zb_4)=0,\]
\begin{eqnarray*}
\lefteqn{d_1d_2(b_1\otimes\alpha\otimes b_2, b_3\otimes b_4)=}\\
&&=d_1(-b_1z\otimes\alpha\otimes b_2+b_1\otimes\alpha\otimes zb_2+\sum\limits_{a\in\bar Q}\epsilon_ab_3a\otimes a^*\otimes b_4\\
&&\quad+\sum\limits_{a\in\bar Q}\epsilon_ab_3\otimes a\otimes a^*b_4,-b_3\mu\otimes b_4+b_1\alpha\otimes b_2-b_1\otimes\alpha b_2)\\
&&=
-b_1z\alpha\otimes b_2+b_1z\otimes\alpha b_2+b_1\alpha\otimes zb_2-b_1\otimes\alpha zb_2+\sum\limits_{a\in\bar Q}\epsilon_ab_3aa^*\otimes b_4\\
&&\quad-\sum\limits_{a\in\bar Q}\epsilon_ab_3a\otimes a^*b_4+\sum\limits_{a\in\bar Q}\epsilon_ab_3a\otimes a^*b_4-\sum\limits_{a\in\bar Q}\epsilon_ab_3\otimes aa^*b_4-b_3z\mu\otimes b_4\\
&&\quad+b_3\mu\otimes zb_4+b_1\alpha z\otimes b_2-b_1\alpha\otimes zb_2-b_1z\otimes \alpha b_2+b_1\otimes z\alpha b_2=0
\end{eqnarray*}
(since $\sum\limits_{a\in\bar Q}\epsilon_aaa^*=z\mu),$
\begin{eqnarray*}
\lefteqn{d_2d_3(b_1\otimes b_2)=}\\
&&=d_2(\sum\limits_{a\in\bar Q}\epsilon_ab_1a\otimes a^*\otimes b_2+\sum\limits_{a\in\bar Q}\epsilon_ab_1\otimes a\otimes a^*b_2,b_1z\otimes b_2-b_1\otimes zb_2)=\\
&&=(-\sum\limits_{a\in\bar Q}\epsilon_ab_1az\otimes a^*\otimes b_2+
\sum\limits_{a\in\bar Q}\epsilon_ab_1a\otimes a^*\otimes zb_2-
\sum\limits_{a\in\bar Q}\epsilon_ab_1z\otimes a\otimes a^*b_2\\
&&\quad+\sum\limits_{a\in\bar Q}\epsilon_ab_1\otimes a\otimes za^*b_2+
\sum\limits_{a\in\bar Q}\epsilon_ab_1za\otimes a^*\otimes b_2+
\sum\limits_{a\in\bar Q}\epsilon_ab_1z\otimes a\otimes a^*b_2\\
&&\quad-\sum\limits_{a\in\bar Q}\epsilon_ab_1a\otimes a^*\otimes zb_2-
\sum\limits_{a\in\bar Q}\epsilon_ab_1\otimes a\otimes a^*zb_2,\\
&&\quad-b_1z\mu\otimes b_2+b_1\otimes z\mu b_2+
\sum\limits_{a\in\bar Q}\epsilon_ab_1aa^*\otimes b_2-
\sum\limits_{a\in\bar Q}\epsilon_ab_1a\otimes a^*b_2\\
&&\quad+\sum\limits_{a\in\bar Q}\epsilon_ab_1a\otimes a^*b_2-
\sum\limits_{a\in\bar Q}\epsilon_ab_1\otimes aa^*b_2)=0,
\end{eqnarray*}
\begin{eqnarray*}
\lefteqn{d_3d_4(b_1\otimes b_2)=d_3(\sum b_1x_i\otimes x_i^*b_2)=}\\
&&=(\sum\limits_{a\in\bar Q}\sum\epsilon_ab_1x_ia\otimes a^*\otimes x_i^*b_2+
\quad\sum\limits_{a\in\bar Q}\sum\epsilon_a b_1x_i\otimes a\otimes a^*x_i^*b_2,\\
&&\quad\sum b_1x_iz\otimes x_i^*b_2-\sum b_1x_i\otimes zx_i^*b_2).
\end{eqnarray*}

Using the trace form, it is easy to show that $\sum x_ia\otimes x_i^*=\sum x_i\otimes ax_i^*$ for any $a\in A$:
\[\sum x_ia\otimes x_i^*=\sum\sum (x_ia,x_j^*)x_j\otimes x_i^*=
\sum\sum x_i\otimes (ax_i^*,x_j)x_j^*=\sum x_i\otimes ax_i^*.\]\\
It follows that $\sum b_1x_iz\otimes x_i^*b_2-\sum b_1x_i\otimes zx_i^*b_2=0$.

Similarly, $\sum x_ia\otimes b\otimes x_i^*=\sum x_i\otimes b\otimes ax_i^*$ for any $a\in A$. Therefore
\[\sum\epsilon_ab_1x_ia\otimes a^*\otimes x_i^*b_2=
\sum\epsilon_ab_1x_i\otimes a^*\otimes ax_i^*b_2=
\sum\underbrace{\epsilon_{a^*}}_{=-\epsilon_a}b_1x_i\otimes a\otimes a^*x_i^*b_2,\] so $d_3d_4=0$.
\begin{eqnarray*}
\lefteqn{d_4d_1[2h](b_1\otimes\alpha\otimes b_2,b_3\otimes b_4)=}\\
&&=d_0(b_1\alpha\otimes b_2-b_1\otimes\alpha b_2+b_3z\otimes b_4-b_3\otimes zb_4)\sum x_i\otimes x_i^*=0.
\end{eqnarray*}

Now we show exactness. Since the complex is periodic, it is enough to show exactness for $C_0$, $C_1$, $C_2$ and $C_3$.

We recall the definition of Anick's resolution \cite{An}. Denote $T_RW$ to be the tensor algebra of a graded $R$-bimodule $W$, $T_R^+W$ its augmentation ideal. Let $L\subset T_R^+W$ be an $R$-graded bimodule and $B=T_RW/(L)$. Then we the following resolution:

\begin{equation}
 B\otimes_RL\otimes_RB\stackrel{\partial}{\rightarrow}B\otimes_RW\otimes_RB\stackrel{f}{\rightarrow}B\otimes_RB\stackrel{m}{\rightarrow}B\rightarrow 0,
\end{equation}

where $m$ is the multiplication map, $f$ is given by
\[
 f(b_1\otimes w\otimes b_2)=b_1w\otimes b_2-b_1\otimes wb_2
\]
and $\partial$ is given by
\[\partial(b_1\otimes l\otimes b_2)=b_1\cdot D(l)\cdot b_2,\]
\begin{eqnarray*}
 D:T_R^+W&\rightarrow&B\otimes_RW\otimes_RB,\\
 w_1\otimes\ldots\otimes w_n&\mapsto&\sum\limits_{p=1}^n(\overline{w_1\otimes\ldots\otimes w_{p-1}})\otimes w_p\otimes(\overline{w_{p+1}\otimes\ldots\otimes w_{n}}),
\end{eqnarray*}
where bar stands for the image in $B$ of the projection map.

In our setting, $W=V\oplus Rz$, $L$ the $R$-bimodule generated by $\sum\limits_{a\in\bar Q}\epsilon_aaa^*-\mu z$ and $\alpha z-z\alpha\,\forall\alpha\in\bar Q$. Then $B=A$.\\
In Anick's resolution, $m=d_0$, $A\otimes_RW\otimes_R A$ can be identified with $C_1$ (via $A\otimes_RA[2]=A\otimes_RRz\otimes_RA$), so that $f$ becomes $d_1$. Then Im$(\partial)$=Im$(d_2)\subset C_1$. This implies exactness in $C_0$ and $C_1$.

For exactness in $2^{nd}$ and $3^{rd}$ term, we show that the complex 
\[
C_4=C_0[2h]\stackrel{d_4}{\rightarrow} C_3\stackrel{d_3}{\rightarrow} C_2\stackrel{d_2}{\rightarrow}C_1\stackrel{d_1}{\rightarrow} C_0\stackrel{d_0}{\rightarrow}A=C_{-1}\rightarrow 0
\]
is selfdual:

By replacing $C_4=C_0[2h]$ by $\bar{C_4}=\text{Im}(d_4)$, we get the complex 
\[
0\rightarrow\bar C_4\stackrel{\bar d_4}{\rightarrow}\bar C_3\stackrel{\bar d_3}{\rightarrow}\bar C_2\stackrel{\bar d_2}{\rightarrow}\bar C_1\stackrel{\bar d_1}{\rightarrow}\bar C_0\stackrel{\bar d_0}{\rightarrow}A\rightarrow 0.
\]

Now, the map $\sum b_1x_i\otimes x_i^*b_2\mapsto b_1b_2$ allows us to identify $Im(d_4)\cong A[2h]$ as $A-$bimodules
so $\bar d_4$ becomes multiplication with $\sum x_i\otimes x_i^*$.

We introduce the following nondegenerate, bilinear forms:\\
On $A\otimes_RA$, let 
\[(x\otimes y,a\otimes b)=Tr(xb)Tr(ya),\] and on $A\otimes_RV\otimes_RA$, we define \[
(x\otimes\alpha\otimes y,a\otimes\beta\otimes b)=Tr(xb)Tr(ya)(\alpha,\beta),\] where we define the form on $V$ by \[(\alpha,\beta)=\epsilon_\beta\delta_{\alpha^*\beta}\]
($\alpha,\beta\in\bar Q$ and $\delta_{ab}=\left\{\begin{array}{cc}1&a=b\\0&a\neq b\end{array}\right.$).

Via the trace form $(x,y)=Tr(xy)$, we can identify $A\cong A^*$, $x\mapsto (x,-)$, and similarly we can use the forms from above to identify $A\otimes_RA\cong (A\otimes_RA)^*$ and $A\otimes_RV\otimes_RA\cong (A\otimes_RV\otimes_RA)^*$, which induces an identification $\bar C_i=\bar C_{3-i}^*$.

We claim the following: $\bar d_0^\star=\bar d_4$, $\bar d_1^\star=-\bar d_3$ and $\bar d_2^\star\iota=\bar d_2$,\\where $\iota(x,y)=(-x,y)$:
\begin{eqnarray*}
(\bar d_4(x),(b_1\otimes b_2))&=&(\sum xx_i\otimes x_i^*,b_1\otimes b_2)=\sum Tr(xx_ib_2)Tr(x_i^*b_1)\\&=&\sum (b_2x,x_i)(x_i^*,b_1)=(b_2x,b_1)=(x,b_1b_2)\\&=&(x,\bar d_0(b_1\otimes b_2)).
\end{eqnarray*}

For $\alpha,\beta\in\bar Q$,
\begin{eqnarray*}
\lefteqn{(-\bar d_3(x\otimes y),(b_1\otimes\alpha\otimes b_2,b_3\otimes b_4))=}\\
&&=((-\sum\limits_{a\in\bar Q}\epsilon_axa\otimes a^*\otimes y-\sum\limits_{a\in\bar Q}\epsilon_ax\otimes a\otimes a^*y,-xz\otimes y+x\otimes zy),\\
&&\quad(b_1\otimes\alpha\otimes b_2,b_3\otimes b_4))\\
&&=-Tr(x\alpha b_2)Tr(yb_1)+Tr(xb_2)Tr(\alpha yb_1)-Tr(xzb_4)Tr(yb_3)\\
&&\quad+Tr(xb_4)Tr(zyb_3)\\
&&=Tr(xb_2)Tr(yb_1\alpha)-Tr(x\alpha b_2)Tr(yb_1)+Tr(xb_4)Tr(yb_3z)\\
&&\quad-Tr(xzb_4)Tr(yb_3)\\
&&=(x\otimes y,b_1\alpha\otimes b_2-b_1\otimes\alpha b_2+b_3z\otimes b_4-b_3\otimes zb_4)\\
&&=(x\otimes y,d_1(b_1\otimes\alpha\otimes b_2,b_3\otimes b_4)),
\end{eqnarray*}
\begin{eqnarray*}
\lefteqn{(\bar d_2(b_1\otimes\alpha\otimes b_2,b_3\otimes b_4),(c_1\otimes\beta\otimes c_2,c_3\otimes c_4))=}\\
&&=((-b_1z\otimes\alpha\otimes b_2+b_1\otimes\alpha\otimes zb_2+\sum\limits_{a\in\bar Q}\epsilon_ab_3a\otimes a^*\otimes b_4\\
&&\quad+\sum\limits_{a\in\bar Q}\epsilon_ab_3\otimes a\otimes a^*b_4,-b_3\mu\otimes b_4+b_1\alpha\otimes b_2-b_1\otimes\alpha b_2),\\
&&\quad(c_1\otimes\beta\otimes c_2,\quad c_3\otimes c_4))\\
&&=-Tr(b_1zc_2)Tr(b_2c_1)(\alpha,\beta)+Tr(b_1c_2)Tr(zb_2c_1)(\alpha,\beta)\\
&&\quad+Tr(b_3\beta c_2)Tr(b_4c_1)-Tr(b_3c_2)Tr(\beta b_4c_1)\\
&&\quad-Tr(b_3\mu c_4)Tr(b_4c_3)+Tr(b_1\alpha c_4)Tr(b_2c_3)-Tr(b_1c_4)Tr(\alpha b_2c_3)\\
&&=Tr(b_1c_2)Tr(b_2c_1z)(\alpha,\beta)-Tr(b_1zc_2)Tr(b_2c_1)(\alpha,\beta)\\
&&\quad-Tr(b_1c_4)Tr(b_2c_3\alpha)+Tr(b_1\alpha c_4)Tr(b_2c_3)\\
&&\quad-Tr(b_3\mu c_4)Tr(b_4c_3)-Tr(b_3c_2)Tr(b_4c_1\beta)+Tr(b_3\beta c_2)Tr(b_4c_1)
\end{eqnarray*}
\begin{eqnarray*}
&&=((b_1\otimes\alpha\otimes b_2,b_3\otimes b_4),(c_1z\otimes\beta\otimes c_2-c_1\otimes\beta\otimes zc_2\\
&&\quad+\sum\limits_{a\in\bar Q}\epsilon_ac_3a\otimes a^*\otimes c_4+\sum\limits_{a\in\bar Q}\epsilon_ac_3\otimes a\otimes a^*c_4,\\
&&\quad-c_3\otimes\mu c_4-c_1\beta\otimes c_2+c_1\otimes\beta c_2))\\
&&=((b_1\otimes\alpha\otimes b_2,b_3\otimes b_4),\bar d_2(-c_1\otimes\beta\otimes c_2,c_3\otimes c_4)).
\end{eqnarray*}

Now, the selfduality of our complex $\bar C_\bullet$ and exactness in $\bar C_0$ and $\bar C_1$ implies exactness in $\bar C_2$ and $\bar C_3$.
\end{proof}
\subsection{Computation of Hochschild cohomology/homology}

Now we use the projective resolution $C_\bullet$ to compute the Hochschild cohomology and homology groups of $A$. Let us write $A^e=A\otimes_RA^{op}$. 
\begin{theorem}\label{cohomology}

The Hochschild cohomology groups of $A$ are:
\begin{eqnarray*}
HH^0(A)&=&Z\textnormal{ (the center of }A),\\
HH^{4n+1}(A)&=&(Z\cap\mu^{-1}[A,A])[-2nh-2],\\
HH^{4n+2}(A)&=&A/([A,A]+\mu Z)[-2nh-2],\\
HH^{4n+3}(A)&=&A_+/[A,A][-2nh-4],\\
HH^{4n+4}(A)&=&Z/A_{top}[-2(n+1)h]
\end{eqnarray*}
where $n\geq0$, and $A_{top}$ is the top-degree part of $A$.
\end{theorem}
\begin{proof}
Apply the functor $Hom_{A^e}(-,A)$ on $C_\bullet$, identify \[Hom_{A^e}(A\otimes_RA,A)\cong A^R\] ($\phi\in Hom_{A^e}(A\otimes_RA,A)$ is determined by
$\phi(1\otimes 1)=a\in A$ and observe $ra=\phi(r\otimes 1)=\phi(1\otimes r)=ar,\,\forall r\in R$. We write $a\circ -$ for $\phi$) and \[
Hom_{A^e}(A\otimes_RV\otimes_RA,A)\cong (A\otimes_RV)^R[-2]
\] 
($\sum\limits_{a\in\bar Q}x_a\otimes a^*$ is identified with the homomorphism $\psi$ which maps each element $1 \otimes a \otimes 1$ to $x_a$ ($a\in\bar Q$), we write $\sum\limits_{a\in\bar Q}x_a\otimes a^*$ for $\psi(-)$)\\
to obtain the Hochschild cohomology complex
\[
\cdots\leftarrow A^R[-2h]\stackrel{d_4^*}{\leftarrow} A^R[-4]\stackrel{d_3^*}{\leftarrow} \begin{array}{ccc}(A\otimes_RV)^R[-4]\\ \oplus\\ A^R[-2]\end{array}\stackrel{d_2^*}{\leftarrow} \begin{array}{ccc} (A\otimes_RV)^R[-2]\\ \oplus \\ A^R[-2]\end{array}\stackrel{d_1^*}{\leftarrow} A^R\leftarrow 0.
\]
\begin{eqnarray*}
\lefteqn{d_1^*(x)(b_1\otimes\alpha\otimes b_2,b_3\otimes b_4)=x\circ d_1(b_1\otimes\alpha\otimes b_2,b_3\otimes b_4)=}\\
&&=b_1\alpha xb_2-b_1x\alpha b_2+b_3zxb_4-b_3xzb_4=b_1[\alpha,x]b_2,
\end{eqnarray*}
so 
\[
 d_1^*(x)=(\sum\limits_{a\in\bar Q}[a,x]\otimes a^*,0).
\]

Let $\alpha=\sum\limits_{a\in\bar Q}r_aa$, $r_a\in R$.
\begin{eqnarray*}
\lefteqn{d_2^*(\sum\limits_{a\in\bar Q}x_a\otimes a^*,0)(b_1\alpha\otimes b_2,b_3\otimes b_4)=(\sum\limits_{a\in\bar Q}x_a\otimes a^*)\circ d_2(b_1\alpha\otimes b_2,b_3\otimes b_4)=}\\
&&=\sum\limits_{a\in\bar Q}(x_a\otimes a^*)\circ (-b_1z\otimes\alpha\otimes b_2+b_1\otimes\alpha\otimes zb_2+\sum\limits_{\beta\in\bar Q}\epsilon_\beta b_3\beta\otimes\beta^*\otimes b_4\\
&&\quad+\sum\limits_{\beta\in\bar Q}\epsilon_\beta b_3\otimes\beta\otimes\beta^*b_4)=\\
&&=\sum\limits_{a\in\bar Q}(-b_1zr_ax_ab_2+b_1r_ax_azb_2)-\sum\limits_{a\in\bar Q}\epsilon_ab_3 a^*x_ab_4+\sum\limits_{a\in\bar Q}\epsilon_ab_3 x_aa^*b_4\\&&=\sum\limits_{a\in\bar Q}\epsilon_ab_3 [x_a,a^*]b_4,
\end{eqnarray*}
so 
\[
d_2^*(\sum\limits_{a\in\bar Q}x_a\otimes a^*,0)=(0,\sum\limits_{a\in\bar Q}\epsilon_a[x_a,a^*]).
\]
\begin{eqnarray*}
\lefteqn{d_2^*(0,y)(b_1\otimes\alpha\otimes b_2,b_3\otimes b_4)=y\circ d_2(b_1\otimes\alpha\otimes b_2, b_3\otimes b_4)=}\\&&=y\circ(-b_3\mu\otimes b_4+b_1\alpha\otimes b_2-b_1\otimes\alpha b_2)=-b_3\mu yb_4+b_1\alpha yb_2-b_1y\alpha b_2, 
\end{eqnarray*}
so
\[
d_2^*(0,y)=(-\sum\limits_{a\in\bar Q}[y,a]\otimes a^*,-\mu y).
\]
Putting this together, we obtain:
\[
d_2^*(\sum\limits_{a\in\bar Q}x_a\otimes a^*,y)=(-\sum\limits_{a\in\bar Q}[y,a]\otimes a^*,-\mu y+\sum\limits_{a\in\bar Q}\epsilon_a[x_a,a^*]).
\]
\begin{eqnarray*}
\lefteqn{d_3^*(\sum\limits_{a\in\bar Q}x_a\otimes a^*,0)(b_1\otimes b_2)=(\sum\limits_{a\in\bar Q}x_a\otimes a^*)\circ d_3(b_1\otimes b_2)=}\\
&&=(\sum\limits_{a\in\bar Q}x_a\otimes a^*)\circ (\sum\limits_{\alpha\in\bar Q}\epsilon_\alpha b_1\alpha\otimes\alpha^*\otimes b_2+\sum\limits_{\alpha\in\bar Q}\epsilon_\alpha b_1\otimes\alpha\otimes\alpha^*b_2,
\\&&\quad b_1z\otimes b_2-b_1\otimes zb_2)=\\
&&=\sum\limits_{a\in\bar Q}(-\epsilon_ab_1a^*x_ab_2+\epsilon_ab_1x_aa^*b_2)=\sum\limits_{a\in\bar Q}\epsilon_ab_1[x_a,a^*]b_2,\\
\lefteqn{d_3^*(0,y)(b_1\otimes b_2)=}\\&&=y\circ (\sum\limits_{\alpha\in\bar Q}\epsilon_\alpha b_1\alpha\otimes\alpha^*\otimes b_2+\sum\limits_{\alpha\in\bar Q}\epsilon_\alpha b_1\otimes\alpha\otimes\alpha^*b_2,b_1z\otimes b_2-b_1\otimes zb_2)
\\&&=b_1zyb_2-b_1yzb_2=0,
\end{eqnarray*}
so we get
\[
d_3^*(\sum\limits_{a\in\bar Q}x_a\otimes a^*,y)=\sum\limits_{a\in\bar Q}\epsilon_a[x_a,a^*].
\]
\[ 
d_4^*(x)(b_1\otimes b_2)=x\circ (\sum b_1x_i\otimes x_i^*b_2)=\sum b_1x_ixx_i^*b_2, 
\]
so
\[
d_4^*(x)=\sum x_ixx_i^*.
\]

Now, we want to compute the Hochschild cohomology (since the complex is periodic, $HH^i(A)=HH^{i+4}(A)[2h]\,\forall i\geq 1$, so it is enough to do the calculations until $HH^4$):

\underline{$HH^0(A)=Z$} (the center of A), since a cocycle $x\in\ker d_1^*$ lies in $A^R$ and has to satisfy $\sum\limits_{a\in\bar Q}[a,x]\otimes a^*=0$, i.e. commute with all $a\in\bar Q$.

\underline{$HH^1(A)=(Z\cap\mu^{-1}[A,A])[-2]:$} ($\sum\limits_{a\in\bar Q}x_a\otimes a^*,y)$ is a cocycle if $\sum\limits_{a\in\bar Q}[y,a]\otimes a^*=0$ (i.e. $y\in Z$) and $y=\mu^{-1}\sum\limits_{a\in\bar Q}\epsilon_a[x_a,a^*]$ (since $\mu$ is invertible) which implies $y\in\mu^{-1}[A,A]$. Since $\sum\limits_{a\in\bar Q}\epsilon_a[a^*,x_a]=0$ implies that $x_a=[a,x]$ (we refer to \cite[Corollary 3.5.]{ELR} where this statement follows from the exactness of the complex in the $1^{st}$ term) for some $x\in A$, and $\sum\limits_{a\in\bar Q}[a,x]\otimes a^*$ lies in Im$d_1^*$, $HH^1(A)$ is controlled only by $y\in (Z\cap\mu^{-1}[A,A])[-2]$. Since $[A(1),A]=[A,A]$, any $y\in (Z\cap\mu^{-1}[A,A])[-2]$ also gives rise to a cocycle.

\underline{$HH^2(A)=A/([A,A]+\mu Z)[-2]:$} An element $(\sum\limits_{a\in\bar Q}x_a\otimes a^*,y)$ is a cocycle if $\sum\limits_{a\in\bar Q}\epsilon_a[x_a,a^*]=0$, so $x_a=[x,a]$ for some $x\in A^R$, (where $x$ is unique up to a central element), so cocycles are of the form $(\sum\limits_{a\in\bar Q}[x,a]\otimes a^*,y)$. The coboundaries are spanned by $(\sum\limits_{a\in\bar Q}[x,a]\otimes a^*,\mu x)$ (where the first component determines $x$ uniquely modulo $Z$) and $(0,\sum\limits_{a\in\bar Q}\epsilon_a[x_a,a^*])$ (where the image is $[A,A]^R$). It follows that 
\[
HH^2(A)=A^R/([A,A]^R+\mu Z)[-2]=A/([A,A]+\mu Z)[-2].
\]

\underline{$HH^3(A)=A_+/[A,A][-4]:$} We denote $A_+$ to be the positive degree part of $A$. $d_4^*(x)=\sum x_ixx_i^*$ is zero if $x$ has positive degree (since $x_ixx_i^*$ exceeds the top degree $2h-4$). 

Observe also that $d_4^*$ injects $R$ into $A_{top}$: \\
Since $A=\oplus e_kAe_j$, we can choose a basis $\{x_i\}$, such that these elements all belong to a certain subspace $e_kAe_j$ for some $k,j$. We denote $\{x_{i'}^{j,k}\}$ the subbasis of $\{x_i\}$ which spans $e_kAe_j$.

Assume that $0=d_4^*(\sum\limits_{j=1}^r\lambda_je_j)$. Then $\forall k$,
\begin{eqnarray*} 
0&=&\sum\limits_{j=1}^r\lambda_j Tr(\sum_i e_kx_ie_jx_i^*)=\sum\limits_{j=1}^r\lambda_j \sum\limits_{i',j,k}\underbrace{(x_{i'}^{j,k},(x_{i'}^{j,k})^*)}_{=1}=\sum\limits_{j=1}^r\lambda_j\dim e_kAe_j\\
&=&\sum\limits_j\lambda_j\sum\limits_d\dim e_kA[d]e_j=H_A(1)_{k,j}=\sum\limits_j\lambda_j\left(\frac{h}{2-C}\right)_{k,j}.
\end{eqnarray*}
The last equality follows from \cite[Theorem 3.2.]{ER}. Since the matrix $\frac{h}{2-C}$ is nondegenerate, all $\lambda_j=0$.

So we see that the images $d_4^*(e_j)$ are nonzero and linearly independent. So the cocycles are the elements in $A_+^R$, and the coboundaries are $\sum\limits_{a\in\bar Q}\epsilon_a[x_a,a^*]$ which generate $[A,A]^R$. Therefore $HH^3(A)=A_+^R/[A,A]^R[-4]=A_+/[A,A][-4].$

\underline{$HH^4(A)=Z/A_{top}[-2h]$:} Since $d_5^*=d_1^*$, the cocycles are the central elements. From the above discussion about the image of $d_4^*$ and the fact that $A_{top}$ is $r$-dimensional, it follows that the coboundaries are the top degree elements of $A$.
\end{proof}

Similarly, we compute the Hochschild homology groups of $A$.
\begin{theorem}
The Hochschild homology groups of $A$ are:
\begin{eqnarray*}
HH_0(A)&=&A/[A,A],\\
HH_{4n+1}(A)&=&A/([A,A]+\mu Z)[2nh+2],\\
HH_{4n+2}(A)&=&(Z\cap\mu^{-1}[A,A])[2nh+2],\\
HH_{4n+3}(A)&=&Z/A_{top}[2nh+4],\\
HH_{4n+4}(A)&=&A_+/[A,A][2(n+1)h].
\end{eqnarray*}
\end{theorem}
\begin{proof}
Apply the functor $(A\otimes_{A^e}- )$ to $C_\bullet$, identify \[A\otimes_{A^e}(A\otimes_RA)\cong A^R\] ($a\otimes (b\otimes c)=cab\otimes 1\otimes1\mapsto cab$ and observe 
\\$\forall a\in A, r\in R: ar=a\otimes (r\otimes 1)=a\otimes (1\otimes r)=ra$) and \[A\otimes_{A^e}(A\otimes_R V\otimes_RA)\cong (A\otimes_RV)^R\] (via $a\otimes(b\otimes\alpha\otimes c)=cab\otimes(1\otimes\alpha\otimes 1)\mapsto cab\otimes\alpha$). 
                        
We get the following periodic complex for computing the Hochschild homology:
\[
\cdots\rightarrow A^R[2h]\stackrel{d_4'}{\rightarrow}A^R[4]\stackrel{d_3'}{\rightarrow}\begin{array}{ccc}(A\otimes_RV)^R[2]\\\oplus\\ A^R[2]\end{array}\stackrel{d_2'}{\rightarrow}\begin{array}{ccc}(A\otimes_RV)^R\\\oplus\\ A^R[2]\end{array}\stackrel{d_1'}{\rightarrow}A^R\rightarrow 0.
\]

The differentials become:
\begin{eqnarray*}
\lefteqn{d_1'(\sum\limits_{a\in\bar Q}x_a\otimes a,y)=1\otimes d_1(\sum\limits_{a\in\bar Q}x_a\otimes a\otimes 1,y\otimes 1)}\\
&&=1\otimes(\sum\limits_{a\in\bar Q}x_aa\otimes 1-\sum\limits_{a\in\bar Q}x_a\otimes a+yz\otimes 1-y\otimes z)=\sum\limits_{a\in\bar Q}[x_a,a],
\end{eqnarray*}
\begin{eqnarray*}
\lefteqn{d_2'(\sum\limits_{a\in\bar Q}x_a\otimes a,y)=1\otimes d_2(\sum\limits_{a\in\bar Q}x_a\otimes a\otimes 1,y\otimes 1)=}\\
&&=1\otimes(-\sum\limits_{a\in\bar Q}x_az\otimes a\otimes 1+\sum\limits_{a\in\bar Q}x_a\otimes a\otimes z+\sum\limits_{\alpha\in\bar Q}\epsilon_\alpha y \alpha\otimes \alpha^*\otimes 1\\
&&\quad+\sum\limits_{\alpha\in\bar Q}\epsilon_\alpha y\otimes \alpha\otimes \alpha^*,
-y\mu\otimes 1+\sum\limits_{a\in\bar Q}x_aa\otimes1-\sum\limits_{a\in\bar Q}x_a\otimes a)\\
&&=(-\sum\limits_{a\in\bar Q}x_az\otimes a+\sum\limits_{a\in\bar Q}zx_a\otimes a+\sum\limits_{\alpha\in\bar Q}\epsilon_\alpha y \alpha\otimes \alpha^*+\sum\limits_{\alpha\in\bar Q}\epsilon_\alpha \alpha^* y\otimes\alpha,\\
&&\quad-y\mu+\sum\limits_{a\in\bar Q}x_aa-\sum\limits_{a\in\bar Q}ax_a)=(\sum\limits_{\alpha\in\bar Q}\epsilon_\alpha[y,\alpha]\otimes\alpha^*,\sum\limits_{a\in\bar Q}[x_a,a]-y\mu),
\end{eqnarray*}
\begin{eqnarray*}
\lefteqn{d_3'(x)=1\otimes d_3(x\otimes 1)=}\\
&&=(1\otimes (\sum\limits_{a\in\bar Q}\epsilon_axa\otimes a^*\otimes 1+\sum\limits_{a\in\bar Q}\epsilon_ax\otimes a\otimes a^*),1\otimes (xz\otimes 1-x\otimes z))\\
&&=(\sum\limits_{a\in\bar Q}\epsilon_axa\otimes a^*+\sum\limits_{a\in\bar Q}\epsilon_aa^*x\otimes a,xz-zx)=(\sum\limits_{a\in\bar Q}\epsilon_a[x,a]\otimes a^*,0),
\end{eqnarray*}
\[
d_4'(x)=1\otimes d_4(x\otimes 1)=1\otimes\sum xx_i\otimes x_i^*=\sum x_i^*xx_i.
\]

Now, we compute the homology (and since the complex is periodic,\\ $HH_i(A)=HH_{i+4}(A)$ for $i>0$, so it is enough to calculate the homology up to $HH_4$):

\underline{$HH_0(A)=A/[A,A]$:} the boundaries are of the form $\sum\limits_{a\in\bar Q}[x_a,a]$, and they generate $[A,A]^R$. So $HH_0(A)=A^R/[A,A]^R=A/[A,A]$ follows.

\underline{$HH_1(A)=A/([A,A]+\mu Z)[2]$:} The cycle condition $\sum\limits_{a\in\bar Q}[x_a,a]=0$ implies $x_a=\epsilon_a[x,a^*]$ for some $x\in A$ (again, we refer to the result $H_1=0$ in\\
\cite[Corollary 3.5.]{ELR}), so the cycles are $(\sum\limits_{a\in\bar Q}\epsilon_a[x,a^*]\otimes a,y)$.\\
The boundaries are of the form $(\sum\limits_{a\in\bar Q}\epsilon_a[x,a^*]\otimes a,\sum\limits_{\alpha\in\bar Q}[x_\alpha,\alpha]+\mu x)$ (where the first component determines $x$ uniquely modulo $Z$. So
\[
HH_1(A)=A^R/([A,A]^R+\mu Z)[2]=A/([A,A]+\mu Z)[2].
\]

\underline{$HH_2(A)=Z\cap\mu^{-1}[A,A][2]$:} 
The cycle conditions are $\sum\limits_{\alpha\in\bar Q}\epsilon_\alpha[y,\alpha]\otimes\alpha^*=0$ (this tells us $y\in Z$) and 
$\sum\limits_{a\in\bar Q}[x_a,a]-y\mu=0$, so $y\in \mu^{-1}[A,A]$ and $x_a$ unique up to an element of the form $\epsilon_a[x,a^*]$ for some $x\in A$. So the cycles are of the form $(\sum\limits_{a\in\bar Q}x_a\otimes  a,y)$, $y\in Z\cap\mu^{-1}[A,A]$, $x_a$ uniquely controlled by $y$ (mod $\epsilon_a[x,a^*]$) , and the boundaries have the form $(\sum\limits_{a\in\bar Q}\epsilon_a[x,a^*]\otimes a,0)$, i.e. homology is controlled only by $y$ now. So $HH_2(A)=Z\cap\mu^{-1} [A,A][2]$.

\underline{$HH_3(A)=Z/A_{top}[4]$:} The cycle condition $\sum\limits_{a\in\bar Q}\epsilon_a[x,a]\otimes a^*=0$ implies that the cycles are the central elements $Z$. The boundaries $\sum x_i^*xx_i$ consist of the top degree part of $A$, so $HH_3(A)=Z/A_{top}[4]$.

\underline{$HH_4(A)=A_+/[A,A][2h]$:} $\ker d_4=A_+^R$, $\textnormal{Im}\, d_5=\textnormal{Im}\, d_1=[A,A]^R$, therefore $HH_4(A)=A_+^R/[A,A]^R=A_+/[A,A]$.
\end{proof}

\subsection{The intersection $Z\cap \mu^{-1}[A,A]$.}
We found $Z\cap\mu^{-1}[A,A]$ as the $(4i+2)-$th homology and $(4i+1)-$th cohomology group, so to understand the (co)homology of $A$ better, we are interested in its structure.

Now, we define the following Hilbert series:
\begin{align*}
q(t)&=h_{Z\cap\mu^{-1}[A,A]}(t),\\
q_*(t)&=h_{A/([A,A]+\mu Z)}(t).
\end{align*}

To relate both to each other, we prove the following
\begin{proposition}\label{pairing}
The trace form defines a nondegenerate pairing \\
$(Z\cap\mu^{-1}[A,A])\times A/([A,A]+\mu Z)\rightarrow\mathbb{C}$.
\end{proposition}
\begin{proof}
Since the trace form is nondegenerate on $A$, it is enough to show that $(Z\cap\mu^{-1}[A,A])^\bot\subset [A,A]+\mu Z$, or equivalently \\
$([A,A]+\mu Z)^\bot\subset Z\cap\mu^{-1}[A,A]$. The latter follows from $[A,A]^\bot\subset Z$, since \[(x,[y_1,y_2])=Tr(x[y_1,y_2])=Tr([x,y_1]y_2)=([x,y_1],y_2)=0\,\forall y_1,y_2\in A\] implies $[x,y_1]=0$, i.e. $x\in Z$.
\end{proof}

\begin{corollary}
$q(t)$ and $q_*(t)$ are palindromes of each other, i.e.\\ $q(t)=t^{2h-4}q_*(1/t)$.
\end{corollary}

Let us define the Hilbert series $p(t)=h_{A/\mu^{-1}[A,A]}(t)$. We recall from \cite[end of section 2.2.]{ELR} that $p(t)=\sum\limits_{i=1}^r(1+t^2+\ldots +t^{2(m_i-1)})$ where the $m_i$ are the exponents of the root system. Since the trace form also defines a nondegenerate pairing $Z\times A/[A,A]\rightarrow\mathbb{C}$ (see \cite[Corollary 2.2.]{ELR}), it follows for the Hilbert series $p_*(t)=h_Z(t)$ that $p(t)=t^{2h-4}p_*(1/t)$. Since $zZ\subset\mu^{-1}[A,A]$ is spanned by even degree elements, we see that $Z$ is generated as a $\mathbb{C}[z]-$module by elements of degree $2(m_i-1)$.
\begin{proposition}
We have
\[
q_*(t)\geq p(t)-\sum\limits_{i=1}^rt^{2(m_i-1)}=\sum\limits_{i=1}^r\left(1+t^2+\ldots +t^{2(m_i-2)}\right).
\]
\end{proposition}
\begin{proof}
From the exact sequence
\[
0\rightarrow Z/(Z\cap\mu^{-1}[A,A])\rightarrow A/\mu^{-1}[A,A]\rightarrow A/(\mu^{-1}[A,A]+Z)\rightarrow0
\]
we obtain the equation \[q_*(t)=p(t)-h_{Z/(Z\cap\mu^{-1}[A,A])}(t).\] Since $zZ\subset\mu^{-1}[A,A]$ ($z=\mu^{-1}\sum\limits_{a\in\bar Q}[a,a^*]\in\mu^{-1}[A,A])$, we have the inequality \[h_{Z/(Z\cap\mu^{-1}[A,A])}(t)\leq h_{Z/zZ}(t)=\sum\limits_{i=1}^rt^{2(m_i-1)},\] and our inequality \[q_*(t)\geq p(t)-\sum\limits_{i=1}^rt^{2(m_i-1)}\] follows.
\end{proof}
\begin{theorem}\label{qp}
The inequality from above is an equality:
\[
q_*(t)=p(t)-\sum\limits_{i=1}^rt^{2(m_i-1)}.
\]
\end{theorem}
We will prove this in the next section where we compute the cyclic homology groups of $A$. From this, we get a result for our intersection space:
\begin{corollary}
$Z\cap\mu^{-1}[A,A]=zZ$.
\end{corollary}
\subsection{Cyclic homology of A}
The Connes differentials $B_i$ (see \cite[2.1.7.]{Lo}) give us an exact sequence
\[
R\stackrel{B_{-1}}{\rightarrow}HH_0(A)\stackrel{B_0}{\rightarrow}HH_1(A)\stackrel{B_1}{\rightarrow}HH_2(A)\stackrel{B_2}{\rightarrow}HH_3(A)\stackrel{B_3}{\rightarrow}HH_4(A)\stackrel{B_4}{\rightarrow}\ldots .
\]
In our case, we have the following exact sequence:
\begin{align*}
R&\stackrel{B_{-1}}{\rightarrow}A/[A,A]\stackrel{B_0}{\rightarrow}A/([A,A]+\mu Z)[2]\stackrel{B_1}{\rightarrow}Z\cap\mu^{-1}[A,A][2]\stackrel{B_2}{\rightarrow}Z/A_{top}[4]\\
&\stackrel{B_3}{\rightarrow}A_+/[A,A][2h]\stackrel{B_4}{\rightarrow}\ldots,
\end{align*}
and the $B_i$ are all degree-preserving. 

We define the \emph{reduced cyclic homology} (see \cite[2.2.13.]{Lo}) 
\begin{align*}
\overline{HC_i}(A)&=\ker (B_{i+1}:HH_{i+1}(A)\rightarrow HH_{i+2}(A))\\&=\textnormal{Im}(B_i:HH_i(A)\rightarrow HH_{i+1}(A)).
\end{align*}

\begin{theorem}
We get the following cyclic homology groups:
\begin{eqnarray*}
\overline{HC_{4n}}(A)&=&A_+/[A,A][2nh],\\
\overline{HC_{4n+1}}(A)&=&0,\\
\overline{HC_{4n+2}}(A)&=&Z/A_{top}[2nh+4],\\
\overline{HC_{4n+3}}(A)&=&0.
\end{eqnarray*}
\end{theorem}
\begin{proof}
First we observe that $B_{4n+3}=0$, since the elements of $Z/A_{top}[4]$ have degree $\leq (2h-6)+4=2h-2$ and the elements in $A_+/[A,A]\,[2h]$ have degree $\geq 2h+1$. So we have for each $n$ the exact sequences
\begin{align*}
0&\rightarrow  \frac{A_+}{[A,A][2nh]}\stackrel{B_{4n}}{\rightarrow}\frac{A}{[A,A]+\mu Z}[2nh+2]\stackrel{B_{4n+1}}{\rightarrow}(Z\cap\mu^{-1}[A,A])[2nh+2]\\
&\stackrel{B_{4n+2}}{\rightarrow}\frac{Z}{A_{top}}[2nh+4]\rightarrow 0.
\end{align*}

The only thing to show is that $W:=\overline{HC_{4n+1}}(A)=\textnormal{Im}B_{4n+1}=0$.
We will use the following theorem from \cite{EG}:
\begin{theorem}
Let $\chi_{\overline{HC}(A)}(t)=\sum a_kt^k$, the Euler characteristic of ${\overline{HC}(A)}$. Then 
\[
\prod\limits_{k=1}^\infty(1-t^k)^{-a_k}=\prod\limits_{s=1}^\infty \det H_A(t^s)=\prod_{s=1}^\infty 
\left(\frac{1-t^{2hs}}{1-t^{2s}}\right)^r\frac{1}{\det (1-Ct^s+t^{2s})},
\]
where $C$ is the adjacency matrix of the quiver $Q$.
\end{theorem}
Since 
\[
\chi_{\overline{HC}(A)}(t)=\frac{1}{1-t^{2h}}(h_{A_+/[A,A]}(t)-h_W(t)+h_{Z/A_{top}}(t)t^4),
\]
to show $W=0$, it is enough to show that if we set 
\[
\frac{1}{1-t^{2h}}(h_{A_+/[A,A]}(t)+h_{Z/A_{top}}(t)t^4)=\sum b_kt^k,
\]
then 
\[
\prod\limits_{k=1}^\infty(1-t^k)^{b_k}=\prod\limits_{s=1}^\infty\left(\frac{1-t^{2s}}{1-t^{2hs}}\right)^r \det(1-Ct^s+t^{2s}).
\]
We have 
\begin{eqnarray*}
h_{A_+/[A,A]}&=&p(t)-r=\sum\limits_{i=1}^r \frac{t^2-t^{2m_i}}{1-t^2} \textnormal{ and} \\
h_{Z/A_{top}}t^4&=&\sum\limits_{i=1}^r\frac{t^{2(m_i-1)}-t^{2h-4}}{1-t^2}t^4=\sum\limits_{i=1}^{r}\frac{t^{2m_i+2}-t^{2h}}{1-t^2}.
\end{eqnarray*}
From these, we get that 
\[
\sum\limits_{k=1}^\infty b_kt^k=(1+t^{2h}+t^{4h}+\ldots)\sum\limits_{i=1}^r (t^2+t^4+\ldots +t^{2m_i-2}+0+t^{2m_i+2}+\ldots+t^{2h-2}),
\]
\begin{align*}
b_k&=0\textnormal{ if }k\textnormal{ is odd}\\
b_{2k}&=\left\{\begin{array}{cc}0&\textnormal{ if }k\textnormal{ is divisible by }h\\r-\#\{i:m_i=p\}&\textnormal{if }k\equiv p\mod h\end{array}\right..
\end{align*}
\[
\prod\limits_{k=1}^\infty(1-t^k)^{b_k}=\prod\limits_{\stackrel{
n}{n\not\equiv 0\mod h
}}
(1-t^{2n})^r/\prod\limits_{\stackrel{n\geq 0}{i\in I}} (1-t^{2(m_i+nh)})
\]

Now, it comes down to showing that 
\[
\prod\limits_{s=1}^{\infty}\det(1-Ct^s+t^{2s})=\prod\limits_{k=1}^{\infty}(1-q^k)^{n_k},
\]
where $q=t^2$ and $n_k=\left\{\begin{array}{cl}0&\textnormal{if }n\textnormal{ is divisible by }h\\-\#\{i:m_i=p\}&\textnormal{if }n\equiv p\mod h\end{array}\right.$ \\(recall that the $m_i$ are the exponents of our root system), for the different Dynkin quivers of type $A_{n-1}$, $D_{n+1}$, $E_6$, $E_7$ and $E_8$. Here we will use the identities for $\det (1-Ct+t^2)=\prod\limits_{j=1}^r(t^2-e^{2\pi im_j/h})$ from \cite[Corollary 4.5.]{RS}.\\
                                          
\underline{Case 1:} $Q=A_{n-1}$\\
The exponents are $1,\ldots ,n-1$ and the Coxeter number is $h=n$.\\
\[
\det (1-Ct+t^2)=\frac{1-t^{2n}}{1-t^2},
\] so if we set
\[
\prod\limits_{k=1}^\infty(1-q^k)^{n_k}=\prod\limits_{s=1}^{\infty}\det(1-Ct^s+t^{2s})=\prod\limits_{s=1}^{\infty}\frac{1-q^{ns}}{1-q^s}
\]
then \[n_k=\left\{\begin{array}{rl}0&\textnormal{if }{n|k}\\-1&\textnormal{if }{n\not| k}\end{array}
\right. .\]\\

\underline{Case 2:} $Q=D_{n+1}$\\
The exponents are $1,3,\ldots ,2n-1,n$ and the Coxeter number is $h=2n$.\\
\[\det (1-Ct+t^2)=\frac{(1-t^4)(1-t^{4n})}{(1-t^2)(1-t^{2n})},\] so 
\[
\prod\limits_{k=1}^\infty(1-q^k)^{n_k}=\prod\limits_{s=1}^{\infty}\det(1-Ct^s+t^{2s})=\prod\limits_{s=1}^{\infty}\frac{(1-q^{2s})(1-q^{2ns})}{(1-q^{s})(1-q^{ns})}
\]
implies that \[n_k=\textnormal{div}(k,2n)-\textnormal{div}(k,n)+\textnormal{div}(k,2)-1,\] where we denote $\textnormal{div}(p,q)=\left\{\begin{array}{rl}1&\textnormal{if }{q|p}\\0&\textnormal{if }{q\not| p}\end{array}
\right.$.\\
\[n_k=\left\{\begin{array}{ll}
0-0+0-1=-1&k\textnormal{ odd}, k \not\equiv0,n\mod 2n\\
0-0+1-1=0 &k\textnormal{ even}, k \not\equiv0,n\mod 2n\\
0-1+1-1=-1&k\textnormal{ even}, k\equiv n\mod 2n\\
0-1+0-1=-2&k\textnormal{ odd}, k\equiv n\mod 2n\\
1-1+1-1=0 &k\equiv 0\mod 2n
\end{array}\right..\]\\

\underline{Case 3:} $Q=E_6$\\
The exponents are $1,4,5,7,8,11$ and the Coxeter number is $h=12$.\\
\[
\det (1-Ct+t^2)=\frac{(1-t^{24})(1-t^4)(1-t^6)}{(1-t^{12})(1-t^8)(1-t^2)},
\] then
\[
\prod\limits_{k=1}^\infty(1-q^k)^{n_k}=\prod\limits_{s=1}^{\infty}\det(1-Ct^s+t^{2s})=\prod\limits_{s=1}^{\infty}\frac{(1-q^{12s})(1-q^{2s})(1-q^{3s})}{(1-q^{6s})(1-q^{4s})(1-q^{s})}
\]
implies \[n_k=\textnormal{div}(k,12)+\textnormal{div}(k,2)+\textnormal{div}(k,3)-\textnormal{div}(k,6)-\textnormal{div}(k,4)-1.\]\\
Observe that if we have a prime factorization $q=a^2b$ ($a,b$ distinct), then \\
\[\textnormal{div}(k,q)+\textnormal{div}(k,a)+\textnormal{div}(k,b)-\textnormal{div}(k,ab)-\textnormal{div}(k,a^2)-1\] is $-1$ if $k$ and $q$ are relatively prime or if $k\equiv la^2\mod q$ ($l\not=0$) and\\$0$ else.

This proves our case for $12=2^2\cdot3$.\\

\underline{Case 4:} $Q=E_7$\\
The exponents are $1,5,7,9,11,13,17$ and the Coxeter number is\\ $h=18$.\\
\[\det (1-Ct+t^2)=\frac{(1-t^{36})(1-t^6)(1-t^3)}{(1-t^{18})(1-t^{12})(1-t^2)},\] so
\[
\prod\limits_{k=1}^\infty(1-q^k)^{n_k}=\prod\limits_{s=1}^{\infty}\det(1-Ct^s+t^{2s})=\prod\limits_{s=1}^{\infty}\frac{(1-q^{18s})(1-q^{3s})(1-q^{2s})}{(1-q^{9s})(1-q^{6s})(1-q^{s})}
\]
implies \[n_k=\textnormal{div}(k,18)+\textnormal{div}(k,3)+\textnormal{div}(k,2)-\textnormal{div}(k,9)-\textnormal{div}(k,6)-1.\]

We use the same argument as above, for $18=2\cdot3^2$.\\

\underline{Case 5:} $Q=E_8$\\
The exponents are $1,7,11,13,17,19,23,29$ and the Coxeter number is\\ $h=30$.\\
\[\det (1-Ct+t^2)=\frac{(1-t^{60})(1-t^{10})(1-t^6)(1-t^4)}{(1-t^{30})(1-t^{20})(1-t^{12})(1-t^2)},\] then
\[
\prod\limits_{k=1}^\infty(1-q^k)^{n_k}=\prod\limits_{s=1}^{\infty}\det(1-Ct^s+t^{2s})=\prod\limits_{s=1}^{\infty}\frac{(1-q^{30s})(1-q^{5s})(1-q^{3s})(1-q^{2s})}{(1-q^{15s})(1-q^{10s})(1-q^{6s})(1-q^{s})}
\]
implies \begin{align*}n_k&=\textnormal{div}(k,30)+\textnormal{div}(k,5)+\textnormal{div}(k,3)+\textnormal{div}(k,2)\\
&\qquad-\textnormal{div}(k,15)-\textnormal{div}(k,10)-\textnormal{div}(k,6)-1.
\end{align*}

We use a similar argument here: If we have a prime factorization $q=abc$ ($a,b,c$ distinct), then
\[\textnormal{div}(k,q)+\textnormal{div}(k,a)+\textnormal{div}(k,b)+\textnormal{div}(k,c)-\textnormal{div}(k,ab)-\textnormal{div}(k,bc)-\textnormal{div}(k,ac)-1\] is $-1$ if $k$ and $q$ are relatively prime and $0$ else.\\

This proves our case for $30=2\cdot3\cdot5$.
\end{proof}

\begin{proof}(of Theorem \ref{qp}):\\
From the isomorphism \[(Z\cap\mu^{-1}[A,A])[2] \stackrel{B_{2}}{\longrightarrow}Z/A_{top}\] we obtain the equation $t^2q(t)=t^4\sum\limits_{i=1}^r(t^{2(m_i-1)}+\ldots+t^{2h-6})$, so \[q(t)=\sum\limits_{i=1}^r(t^{2m_i}+\ldots+t^{2h-4}).\] Recall the duality of exponents, i.e. $m_{r+1-i}=h-m_i$. Then we get
\begin{eqnarray*}
q_*(t)&=&t^{2h-4}q(1/t)=t^{2h-4}\sum\limits_{i=1}^r(t^{-2m_i}+\ldots+t^{-2h+4})\\
&=&t^{2h-4}\sum\limits_{i=1}^r(t^{-2(h-m_i)}+\ldots+t^{-2h+4})\\
&=&\sum\limits_{i=1}^r(1+\ldots+t^{2(m_i-2)})=p(t)-\sum\limits_{i=1}^rt^{2(m_i-1)}.
\end{eqnarray*}
\end{proof}
\section{Universal deformation of $A$}
\begin{definition}
For any weight $\lambda=(\lambda_i)$, we define the algebra
\[
A(\lambda)=P[z]/\left(\sum\limits_{a\in\bar Q}[a,a^*]=z\mu+\sum\limits_{i=1}^r\lambda_ie_i\right)
\]
and introduce a deformation parametrized by formal variables $c_i^j$, $1\leq i\leq r$, $1\leq j\leq h-1$:
\[A(\lambda)_c=P[z][[c]]/\left(\sum\limits_{a\in\bar  Q}[a,a^*]=z\mu+\sum\limits_{i=1}^r\lambda_ie_i+\sum\limits_{i=1}^r\sum\limits_{j=1}^{h-1}c_i^jz^je_i\right).\]
\end{definition}
\begin{theorem}
This deformation is flat $\forall\lambda\in R$, i.e. $A(\lambda)_c$ is free over $\mathbb{C}[[c]]$, and
\[
A(\lambda)_c/(c)=A(\lambda).
\]
\end{theorem}
\begin{proof}
It is sufficient to check flatness for generic $\lambda$. From \cite[end of section 3.2.]{ER}, we know that for generic $\lambda$, $A(\lambda)=\oplus \text{End} V_\alpha$ is a semisimple algebra. So it suffices to show that the representation $V_\alpha$ can be deformed to a representation of $A(\lambda)_c$ for all $\lambda$. 

We recall from \cite[Theorem 4.3.]{CBH} that $\forall\beta\in R$, such that $\beta\cdot\alpha=0$, it exists an $\alpha-$dimensional irreducible representation $V_\alpha$ of $P$, such that \[\sum\limits_{a\in\bar Q}[a,a^*]=\sum\limits_{i=1}^r\beta_ie_i.\]

If we set $z=\gamma\in\mathbb{C}$ in $A(\lambda)_c$, then the relation becomes 
\[
\sum\limits_{a\in\bar Q}[a,a^*]=\sum\limits_{i=1}^re_i(\lambda_i+\gamma(\mu_i+c_i^1)+\gamma^2 c_i^2+\ldots).
\]
Then for $\alpha=\sum\limits_{i=1}^r \alpha_i\epsilon_i$, since the trace of $[a,a^*]$ is zero, the condition to have an $\alpha-$dimensional representation of $A(\lambda)_c$ (i.e. a representation of $P$ satisfying the above relation) is \[\sum\limits_{i=1}^r\alpha_i(\lambda_i+\gamma(\mu_i+c_i^1)+\gamma^2 c_i^2+\ldots)=0.\] By Hensel's lemma, this equation in $\mathbb{C}[[c]]$ has a unique solution $\gamma$, such that its constant term $\gamma_0\in\mathbb{C}$ satisfies $\sum\limits_{i=1}^r\alpha_i(\lambda+\gamma_0)=0\Rightarrow\gamma_0=-\frac{\sum \alpha_i\lambda_i}{\sum \alpha_i}$.
\end{proof}
In particular, if we treat $\lambda$ as formal parameter, then $A(\lambda)_c$ is a flat deformation of $A(0)$.

Let $E$ be the linear span of $z^je_i$, $0\leq j\leq h-2$, $1\leq i\leq r$. From \cite[Proposition 2.4.]{ELR} we know that the projection map $E\rightarrow A/[A,A]$ is surjective. Then the deformation $A(\lambda)_c$ is parametrized by $E$ which gives us a natural map $\eta:E\rightarrow HH^2(A)$. On the other hand, the isomorphism $HH^2(A)=A/([A,A]+\mu Z)$ in Theorem \ref{cohomology} induces a projection map $\theta:E\rightarrow HH^2(A)$.
\begin{proposition}\label{HH^2}
The maps $\theta,\eta:E\rightarrow HH^2(A)$ are identical.
\end{proposition}
\begin{proof}
We have the following commutative diagram which connects our periodic projective resolution with the bar resolution of $A$,

$$
\CD
\begin{array}{ccc}
(A\otimes V\otimes A[2])\\
\oplus\\
(A\otimes A[2])
\end{array}
@> d_2>> 
\begin{array}{ccc}
(A\otimes V\otimes A)\\
\oplus\\
(A\otimes A[2])
\end{array}
@>d_1>>A\otimes A@>d_0>>A\\
@Vf_2VV@VVf_1V@|@|\\
A^{\otimes 4}@>>\tilde d_2>A^{\otimes 3}@>>\tilde d_1>A^{\otimes 2}@>>\tilde d_0>A,
\endCD
$$

where we define 
\[
f_1(b_1\otimes\alpha\otimes b_2,b_3\otimes b_4)=b_1\otimes\alpha\otimes b_2+b_3\otimes z\otimes b_4
\]
and
\[
 f_2(b_1\otimes\alpha\otimes b_2,b_3\otimes b_4)=-b_1\otimes z\otimes\alpha\otimes b_2+b_1\otimes\alpha\otimes z\otimes b_2+\sum\limits_{a\in Q}\epsilon_ab_3\otimes a\otimes a^*\otimes b_4.
\]

Let us check the commutativity of the diagram:
\begin{eqnarray*}
\lefteqn{\tilde d_1f_1(b_1\otimes\alpha\otimes b_2,b_3\otimes b_4)=}\\&&=b_1\alpha\otimes b_2-b_1\otimes\alpha b_2+b_3 z\otimes b_4-b_3\otimes zb_4=d_1(b_1\otimes\alpha\otimes b_2,b_3\otimes b_4),
\end{eqnarray*}
\begin{eqnarray*}
\lefteqn{f_1d_2(b_1\otimes\alpha\otimes b_2,b_3\otimes b_4)=}\\
&&=f_1(-b_1z\otimes\alpha\otimes b_2+b_1\otimes\alpha\otimes zb_2+\sum\limits_{a\in\bar Q}\epsilon_ab_3a\otimes a^*\otimes b_4\\
&&\quad+\sum\limits_{a\in\bar Q}\epsilon_ab_3\otimes a\otimes a^*b_4,
-b_3\mu\otimes b_4+b_1\alpha\otimes b_2-b_1\otimes \alpha b_2)\\
&&=-b_1z\otimes\alpha\otimes b_2-b_1\otimes z\otimes\alpha b_2+b_1\alpha\otimes z\otimes b_2+b_1\otimes\alpha\otimes zb_2\\
&&\quad+\sum\limits_{a\in\bar Q}\epsilon_ab_3a\otimes a^*\otimes b_4-b_3\otimes z\mu\otimes b_4+\sum\limits_{a\in\bar Q}\epsilon_ab_3\otimes a\otimes a^*b_4\\
&&=\tilde d_2f_2(b_1\otimes\alpha\otimes b_2,b_3\otimes b_4).
\end{eqnarray*}

We apply $Hom_{A^e}(-,A)$ to the above diagram:
$$
\CD
\begin{array}{ccc}
(A\otimes V)^R[2]\\
\oplus\\
(A^R[2])
\end{array}
@< d_2^* <<
\begin{array}{ccc}
(A\otimes V)^R[2]\\
\oplus\\
A^R
\end{array}
@<d_1^*<< A^R \\
@Af_2^*AA @AAf_1^*A @|\\
Hom_{A^e}(A^{\otimes 4},A) @<<(\tilde d_2)^*< Hom_{A^e}(A^{\otimes 3},A) @<< (\tilde d_1)^*< A^{R}.
\endCD
$$

The map $f_2^*$ induces a natural isomorphism on $HH^2(A)$, so via this identification we want to prove that $f_2^*\eta=\theta$.\\

The element $\gamma:=\sum\gamma_i^j z^je_i$, $\gamma_i^j\in\mathbb{C}$ defines the 1-parameter deformation \[A^\gamma=A[[\hbar]]/\sum\limits_{a\in Q}[a,a^*]=z+\hbar(\sum\limits_{i=1}^r\gamma_i^0e_i+\sum\limits_{i=1}^r\sum\limits_{j=1}^{h-2}\gamma_i^jz^je_i),\] so the cocycle $\eta(\gamma)$ is defined to be a bilinear map $m$ on $A\times A$ (where we identify \[Hom_{A^e}(A^{\otimes 4},A)=Hom_k(A\otimes A,A)\] here), such that for $a,b\in A$,
\[
a\ast b\equiv ab+\hbar m(a,b) \mod \hbar^2
\]
where $''\ast''$ is the product in $A^\gamma$.
This gives us:
\begin{eqnarray*}
\lefteqn{f_2^*\eta(\gamma)(b_1\otimes\alpha\otimes b_2,b_3\otimes b_4)=\eta(\gamma) f_2(b_1\otimes\alpha\otimes b_2,b_3\otimes b_4)=}\\
&&=\eta(\gamma)(-b_1\otimes z\otimes\alpha\otimes b_2+b_1\otimes\alpha\otimes z\otimes b_2+\sum\limits_{a\in Q}\epsilon_ab_3\otimes a\otimes a^*\otimes b_4)\\
&&=b_1(m(z,\alpha)-m(\alpha,z))b_2+b_3(\sum\limits_{a\in Q}m(a,a^*)-m(a^*,a))b_4\\
&&=b_3(\sum\limits_{i=1}^r\gamma_i^0e_i+\sum\limits_{i=1}^r\sum\limits_{j=1}^{h-2}\gamma_i^jz^je_i)b_4=\theta(\gamma)(b_1\otimes\alpha\otimes b_2,b_3\otimes b_4).
\end{eqnarray*}

We obtain the second to last equality by:\\
$0+\hbar(m(z,\alpha)-m(\alpha,z))=z\ast \alpha-\alpha\ast z=0$ and\\
\begin{align*}
z\mu+\hbar(\sum\limits_{a\in Q}m(a,a^*)-m(a^*,a))&=\sum\limits_{a\in Q}(a\ast a^*-a^*\ast a)\\
&=z\mu+\hbar\left(\sum\limits_{i=1}^r\gamma_i^0e_i+\sum\limits_{i=1}^r\sum\limits_{j=1}^{h-2}\gamma_i^jz^je_i\right).
\end{align*}
This finishes our proof that $f_2^*\eta=\theta$.
\end{proof}

We see that the map $E\rightarrow HH^2(A)$ induced by the deformation $A(\lambda)_c$ is just the projection map. From this we can derive the universal deformation of $A$ very easily.

Let $E^{'}\subset E$ be the subspace which is complimentary to \\$\ker(\theta:E\rightarrow A/([A,A]+\mu Z))$ with basis $w_i,\ldots,w_s$, and choose formal parameters $t_i,\ldots,t_s$. The subdeformation $A'$ of $A$, parametrized by $E'\subset E$ is: 
\[
A^{'}=P[z][[t_1,\ldots,t_{s}]]/\left(\sum\limits_{a\in\bar Q}[a,a^*]=\mu z+\sum\limits_{i=1}^{s}t_iw_i\right).
\]
\begin{theorem}
$A^{'}$ is the universal deformation of $A$.
\end{theorem}
\begin{proof}
$\eta:E'\rightarrow HH^2(A)$ is the map induced by the deformation $A'$. Since $\theta$ induces an isomorphism $E'\stackrel{~}{\rightarrow}A/([A,A]+\mu Z)=HH^2(A)$, by Proposition \ref{HH^2} $\eta$ is an isomorphism and therefore induces a universal deformation.
\end{proof}


\begin{thebibliography}{999}

\bibitem[An]{An} D. J. Anick, \emph{On the homology of associative algebras} Trans. Amer. Math. Soc. \textbf{296} (1986), 641-659

\bibitem[CBH]{CBH} 
W. Crawley-Boevey and M. P. Holland: \emph{Noncommutative deformations of Kleinian singularities}, Duke Mathematical Journal, Vol. 92, No. 3 (1998)

\bibitem[EG]{EG} P. Etingof and V. Ginzburg: \emph{Noncommutative complete intersections and matrix integrals}, math/0603272

\bibitem[ER]{ER} P. Etingof and E. Rains: \emph{Central extensions of preprojective algebras, the quantum Heisenberg algebra, and 2-dimensional complex reflection groups}, math/0503393

\bibitem[ELR]{ELR} P. Etingof, F. Latour and E. Rains: \emph{On central extensions of preprojective algebras}, math/0606403

\bibitem[ES]{ES} K. Erdmann, N. Snashall,
{\em  On Hochschild cohomology of preprojective algebras. I, II.}
  J. Algebra  \textbf{205}  (1998),   391--412, 413--434.

\bibitem[Lo]{Lo} J.-L. Loday: \emph{Cyclic Homology}, Grundreihen der mathematischen Wissenschaften, Vol. 301; A Series of Comprehensive Studies in Mathematics, Springer-Verlag (1992)

\bibitem[RS]{RS} R. Suter: \emph{Quantum affine Cartan matrices, Poincare series of binary polyhedral groups and reflection representations}, math/0503542
\end{thebibliography}
\end{document}